\documentclass[twoside, 10pt]{article}
\usepackage{amssymb, amsmath, mathrsfs}
\usepackage{graphicx}
\usepackage{color}

\input{epsf}

\newcommand{\bd}{\begin{description}}
\newcommand{\ed}{\end{description}}
\newcommand{\be}{\begin{enumerate}}
\newcommand{\ee}{\end{enumerate}}
\newcommand{\bi}{\begin{itemize}}
\newcommand{\ei}{\end{itemize}}
\newcommand{\bt}{\begin{tabular}}
\newcommand{\et}{\end{tabular}}
\newcommand{\beq}{\begin{equation}}
\newcommand{\eeq}{\end{equation}}
\newcommand{\beqs}{\begin{eqnarray*}}
\newcommand{\eeqs}{\end{eqnarray*}}
\newcommand{\proof}{\noindent \textit{Proof: ~}}

\newcommand{\qed}{{\hfill \bf $\Box$}\\ \vspace{0.1in}}

\newcommand{\flr}[1]{\left\lfloor #1 \right\rfloor}
\newcommand{\ceil}[1]{\left\lceil #1 \right\rceil}

\newcommand{\g}{{\mathfrak g}}

\newcommand{\ind}{\text{ind }}

\newtheorem{theorem}{Theorem} 
\newtheorem{lemma}[theorem]{Lemma} 
\newtheorem{definition}[theorem]{Definition}

\newtheorem{fact}[theorem]{Fact}

\newtheorem{conjecture}[theorem]{Conjecture}

\begin{document}
\title{\bf Families of Frobenius seaweed Lie algebras}
\author{Vincent Coll\footnote{Department of Mathematics, Lehigh University, Bethlehem, PA, USA}, Colton Magnant\footnote{Department of Mathematical Sciences, Georgia Southern University, Statesboro, GA, USA.}, Hua Wang\footnotemark[2]}

\maketitle
 \thispagestyle{empty}
 \pagestyle{myheadings}
 \markboth{V. Coll, C. Magnant, H. Wang}{Families of Frobenius seaweed Lie algebras}

\begin{abstract}
We extend the set of known infinite families of Frobenius seaweed Lie subalgebras of $\mathfrak{sl}_{n}$ to include a family which is the first non-trivial general family containing algebras whose associated meanders have an arbitrarily large number of parts.
\end{abstract}

{\bf Keywords:} Biparabolic, Frobenius, Lie algebra, Meander

\section{Introduction}

Meanders were introduced by Dergachev \& Kirillov in \cite{DK} as planar graph representations of seaweed subalgebras of $\mathfrak{sl}_{n}$. In \cite{DK}, the authors provided a combinatorial method of computing the index of such algebras from the number and type of the connected components in their associated meander graphs.  Of particular interest are those seaweed algebras whose meander graphs consist of a single path and so, by the algorithm of Dergachev \& Kirillov, must have index $0$.  More generally, algebras with index $0$ are called \emph{Frobenius}.  Such algebras are of critical importance to deformation and quantum group theory because of their relation to the classical Yang-Baxter equation (CYBE). Suppose $B_F(-,-)$ is non-degenerate and let $M$ be the matrix of $B_F(-.-)$ relative to some basis $\{x_1,\ldots , x_n\}$ of $\mathfrak{g}$. Belavin \& Drinfel'd showed that $r=\sum_{i,j}(M^{-1})_{ij}x_i\wedge x_j$ is a (constant) solution of the CYBE, see \cite{BD}. Thus, each pair consisting of a Lie algebra $\mathfrak{g}$ together with functional $F\in \mathfrak{g}^*$ such that $B_F$ is non-degenerate provides a solution to the CYBE.  See \cite{GG:Boundary} and \cite{GG:Frob} for examples.

The index of a semisimple Lie algebra $\g$ is equal to its rank.  Therefore, such algebras can never be Frobenius. However, there always exist Frobenius subalgebras of $\g$. In particular, many amongst the class of seaweed subalgebras of $\mathfrak{sl}_{n}$, such as the Borel and Cartan subalgebras, are Frobenius.

Combinatorial methods have been developed to compute the index but unfortunately these methods are difficult to apply in practice \cite{Joseph, Khoroshkin, Ooms, Panyushev, Stolin, Tauvel-Yu}.  This may explain why so few explicit examples of nontrivial families are known.  Following Elashvili \cite{Elashvili2} and Coll et al. \cite{CGM11} (see Theorems~\ref{Thm:3-Blocks} and~\ref{Thm:4-Blocks}) we present a simple number-theoretic description of a new infinite family.  This family is the first non-trivial general family containing algebras whose associated meanders have an arbitrarily large number of parts.  

\section{Definitions}

Let $\mathfrak{g}$ be a Lie algebra over a field of characteristic zero.  For any functional $f \in \mathfrak{g}^*$, one can associate the \emph{Kirillov form} $B_f (x,y) = f[x,y]$ which is skew-symmetric and bilinear.  The \emph{index} of $\mathfrak{g}$ is defined to be 
$$
\ind \mathfrak{g} = \min_{f\in \mathfrak{g}^*}\dim (\ker (B_f)).
$$ 
The Lie algebra $\mathfrak{g}$ is {\it Frobenius} if $\ind \mathfrak{g} =0$ (equivalently if $B_{f}$ is nondegenerate).  We consider a subalgebra of seaweed type, defined by Dergachev \& Kirillov in \cite{DK} as follows.

\begin{definition}
Let ${\bf k}$ be an arbitrary field of characteristic $0$ and $n$ a positive integer.  Fix two ordered partitions $\{ a_{i} \}^{k}_{i = 1}$ and $\{ b_{j} \}^{\ell}_{j = 1}$ of the number $n$.  Let $\{ e_{i} \}^{n}_{i = 1}$ be the standard basis in ${\bf k}^{n}$.  A subalgebra of $\rm{M}_{n \times n}({\bf k})$ that preserves the vector spaces $\{V_{i} = \rm{span}(e_{1}, \dots, e_{a_{1} + \dots + a_{i}}) \}$ and $\{W_{j} = \rm{span}(e_{b_{1} + \dots + b_{j} + 1}, \dots, e_{n}) \}$ is called a subalgebra of \emph{seaweed} type due to the suggestive shape of the subalgebra in the total matrix algebra.  See Figure~\ref{MatrixFig}.
\end{definition}

\begin{figure}[ht!]
\begin{center}
 \epsfbox{Matrix.1}
 \caption{Matrix for $a_{1} = 1$, $a_{2} = 2$, $a_{3} = 3$, $b_{1} = 2$ and $b_{2} = 4$. \label{MatrixFig}}
\end{center}
\end{figure}

To each seaweed algebra, a planar graph, called a \emph{meander}, can be constructed as follows: Line up $n$ vertices and label them with the integers $1, 2, \dots, n$ and partition this set into two ordered partitions, called \emph{top} and \emph{bottom}.  For each part in the top (likewise bottom) we build up the graph by adding edges in the same way.  This involves adding the edge from the first vertex of a part to the last vertex of the same part drawn concave down (respectively concave up in the bottom part case).  The edge addition is then repeated between the second vertex and the second to last and so on within each part of both partitions.

\begin{figure}[ht!]
\begin{center}
 \epsfbox{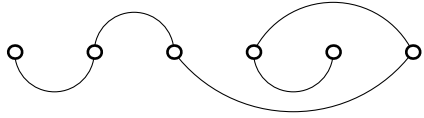}
 \caption{Meander of type $\frac{1|2|3}{2|4}$. \label{MatrixMeanderFig}}
\end{center}
\end{figure}

The parts of these partitions are called \emph{blocks}.  With top blocks of sizes $a_{1}, a_{2}, \dots, a_{\ell}$ and bottom blocks of sizes $b_{1}, b_{2}, \dots, b_{m}$, we say such a meander has \emph{type} $\frac{a_{1}|a_{2}|\dots|a_{\ell}}{b_{1}|b_{2}| \dots |b_{m}}$.  The meander corresponding to the matrix in Figure~\ref{MatrixFig} can be seen in Figure~\ref{MatrixMeanderFig}.  All meanders discussed in this work have multiple top blocks and only one bottom block.  Thus, a meander of type $\frac{a_{1}|a_{2}|\dots|a_{n}}{b}$ will be written simply as $a_{1}|a_{2}|\dots|a_{n}$.  

\noindent \textit{Remark:}  A subalgebra $\mathfrak{p}$ of a semisimple algebra $\mathfrak{g}$ is called \emph{parabolic} if it contains a maximum solvable subalgebra (a Borel subalgebra) of $\mathfrak{g}$.  From the perspective of meanders, a subalgebra is parabolic precisely when the meander has only one block on the bottom.  A subalgebra is \emph{maximal parabolic} if the meander has exactly one bottom block and two top blocks.  Joseph, in \cite{Joseph}, refers to seaweed algebras as \emph{biparabolic} since they are the intersection of two parabolic subalgebras whose sum is $\mathfrak{g}$.

\section{Preliminary results}



The following theorem enables us to easily determine whether a meander with three blocks is Frobenius.  This result was originally established by Elashvili in \cite{Elashvili2}.  Unfortunately, that work remains unpublished.  The theorem was reproven by Coll et al. in \cite{CGM11} using very different techniques.

\begin{theorem}[Elashvili \cite{Elashvili2}, Coll et al. \cite{CGM11}]\label{Thm:3-Blocks}
The meander of type $a|b$ is Frobenius if and only if $\gcd(a, b) = 1$.
\end{theorem}

It is also known precisely when a meander containing four blocks is Frobenius.  This result will be used in our proofs.

\begin{theorem}[Coll et al. \cite{CGM11}]\label{Thm:4-Blocks} 
The meander of type $a|b|c$ is Frobenius if and only if $\gcd(a+b,b+c)=1$.
\end{theorem}

Since the standard drawing of a meander is a planar embedding of the meander graph, a cycle in a meander is a closed curve in the plane.  We recall the classical Jordan Curve Theorem which will help in the proof of our main result.

\begin{theorem}[Jordan Curve Theorem]\label{theorem:Jordan}
Every cycle in a planar embedding of a graph separates the plane into two distinct regions, an interior and an exterior.
\end{theorem}

Our first lemma provides a method of reduction in the size of the blocks.

\begin{lemma}\label{lemma:reduction}
For an even integer $a$ and an odd integer $b$, a meander with $k$ blocks $a|a|\dots |a|b$ has index $m$ if and only if the meander $a|a|\dots |a|(b + 2a)$ with the same number of blocks also has index $m$.
\end{lemma}

\proof
The proof of this lemma consists of simply showing that any path and cycle structure is preserved as the meander is transformed between the two meanders $M$ and $M'$ described below. 
Let $A_{1}, A_{2}, \dots, A_{k - 1}$ denote the $k - 1$ blocks of order $a$ (appearing in this order) with $A_{1} = u_{1}, u_{2}, \dots, u_{a}$ and let $B^{+}$ denote a block of order $b + 2a$ with vertices $v_{1}, v_{2}, \dots, v_{b + 2a}$.  Then all bottom edges of $A_{1}$ go to the rightmost $a$ vertices of $B^{+}$, namely the vertices $v_{b + a + 1}, v_{b + a + 2}, \dots, v_{b + 2a}$.  Using the top edges, these vertices are adjacent to the leftmost $a$ vertices of $B^{+}$, namely the vertices $v_{1}, v_{2}, \dots, v_{a}$.

\begin{figure}[ht!]
\begin{center}
 \epsfxsize 4in
 \epsfbox{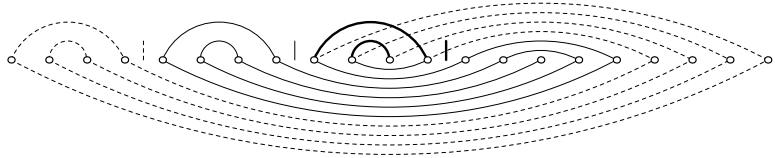}
 \caption{Meanders $M$ (dotted) and $M'$ (thick). \label{Figure:reduction}}
\end{center}
\end{figure}

For each $1 \leq i \leq \frac{a}{2}$, we may replace the path $v_{i}v_{b + 2a - i + 1}u_{i}u_{a - i + 1}v_{b + a + i}v_{a - i + 1}$ with a new edge $v_{i}v_{a - i + 1}$ where all internal vertices of the path are removed from the meander.  Define $A_{1}'$ to be the vertices $v_{1}, v_{2}, \dots, v_{a}$ and define $B$ to be the vertices $v_{a + 1}, v_{a + 2}, \dots, v_{b + a}$.  Note that all the vertices of $A_{1}$ along with $v_{b + a + 1}, v_{b + a + 2}, \dots, v_{b + 2a}$ have been removed.  Call this new graph $M'$.  Then $M'$ is again a meander (see Figure~\ref{Figure:reduction} where the thick edges are only in $M'$ and the dotted edges are only in $M$) and, since we have replaced paths with single edges, all path and cycle structure has been preserved.  Thus, the index is preserved.  This process can easily be reversed to produce $M$ from $M'$.  \qed

In particular, Lemma~\ref{lemma:reduction} means that, when considering whether a general meander of type $a|a|\dots |a|b$ is Frobenius, we may assume $b < 2a$.

Following from the same proof, we also state the following more general version of Lemma~\ref{lemma:reduction}.

\begin{lemma}\label{GenReduction}
For any meander $M$ of type $a_{1}|a_{2}|\dots |a_{t}$ where $a_{t} \geq 2a_{1}$, $M$ has the same index as the meander $a_{2}|a_{3}|\dots |a_{t - 1}|a_{1}|(a_{t} - 2a_{1})$ where the last term is ignored if $a_{t} - 2a_{1} = 0$.
\end{lemma}

\begin{figure}[ht!]
\begin{center}
 \epsfxsize 4in
 \epsfbox{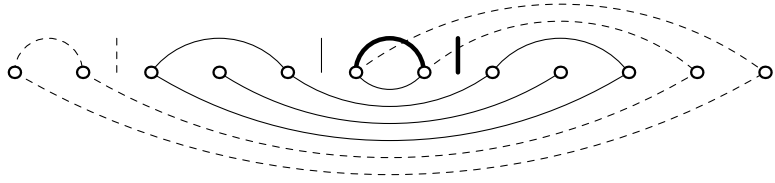}
 \caption{Contracting $2|3|7$ to $3|2|3$. \label{Figure:reduction2}}
\end{center}
\end{figure}

In Figure~\ref{Figure:reduction2}, we see an application of Lemma~\ref{GenReduction} in a meander that is clearly not Frobenius.

\section{A new family of Frobenius meanders}

In the spirit of Theorems~\ref{Thm:3-Blocks} and~\ref{Thm:4-Blocks}, we seek relatively prime conditions on the sizes of the blocks in order to determine whether a given meander is Frobenius.

Call a segment between consecutive vertices in the drawing of a meander an \emph{end-segment} if there is an edge of the meander connecting the two vertices on either end of this segment.  Call a segment a \emph{top-end-segment} if the segment is an end-segment and the corresponding edge is a top-edge.  A segment gets \emph{mapped} by the meander by following either the bottom edges or the top edges on either side of the segment.

\begin{theorem}
A meander of type $a|a|\dots |a|b$ where $a$ is even and $\gcd(a, b) = 1$ is Frobenius.
\end{theorem}

\proof
By Theorem~\ref{Thm:4-Blocks}, we may assume there are at least $3$ blocks of size $a$ in this meander.

\begin{figure}[ht!]
\begin{center}
 \epsfxsize 4in
 \epsfbox{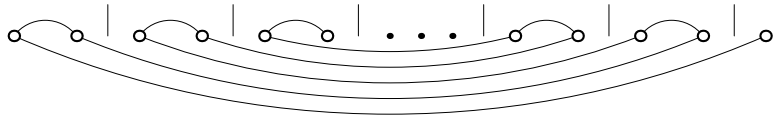}
 \caption{A meander of type $2|2| \dots |2|1$. \label{b-is-1}}
\end{center}
\end{figure}

\begin{figure}[ht!]
\begin{center}
 \epsfxsize 4in
 \epsfbox{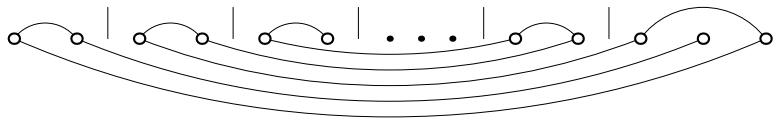}
 \caption{A meander of type $2|2| \dots |2|3$. \label{b-is-3}}
\end{center}
\end{figure}

First suppose $a = 2$.  Then by Lemma~\ref{lemma:reduction}, we know $b \leq 3$.  Then it is easy to see from Figures~\ref{b-is-1} and~\ref{b-is-3} that this meander is Frobenius.

Thus, we may assume $a \geq 4$.  If there was a cycle, Theorem~\ref{theorem:Jordan} implies that there exists a segment between vertices that does not map (following edges of the meander) to the exterior face.  Since we show that this is not the case, there must not be a cycle in the meander.  This proof involves considering the segments between vertices in the meander graph and showing that each top-end-segment must be mapped by the meander to the exterior face.  We then consider any segment and show that it either maps to a top-end-segment or to the exterior face.

Intuitively, the proof makes use of the following heuristic: We visualize the meander as an object (much like a shell) with openings in the top between the blocks.  If water is poured into these openings, we ask whether the water permeates all areas inside the shell.  If so, the meander is Frobenius and if not, it contains a cycle and is therefore not Frobenius.

Label the segments between consecutive vertices of the meander in the natural (drawn) order from $1$ up to $ka + b - 1$ where $k$ is the number of even blocks.  Let $c = \frac{a}{2}$ and note that, since $c \geq 2$, $c$ also shares no common factor with $b$.  Then any top-end-segment must be labeled with $\omega c$ where $\omega$ is an odd positive integer.  Furthermore, the exterior face is accessible via any segment labeled $ia = 2ic$ for any positive integer $i \leq k$.

For any segment labeled $x$, following the bottom mapping by the meander yields $2kc + b - x$ since $c = \frac{a}{2}$.  We denote this by 
$$
x \overset{b}{\longrightarrow} 2kc + b - x.
$$
When following the top mapping by the meander, we have two cases.  If we are in the odd group (of size $b$), we would map $x$ to $4kc + b - x$ so we denote this by 
$$
x \overset{ot}{\longrightarrow} 4kc + b - x.
$$
Suppose now we map within an even group (of size $a$).  If we start in the $\ell^{th}$ such group, we would map $x$ to $(2\ell - 1)2c - x$ which is denoted by
$$
x \overset{et}{\longrightarrow} (2\ell - 1)2c - x.
$$
These maps will be called \emph{arrow maps}.

Our first goal is to show that, starting with a top-end-segment labeled $\omega c$ and following the mapping by the meander, we would reach the exterior face (represented by a segment labeled $2ic$ for some positive integer $i$ as above) before ever reaching another top-end-segment.  In order to accomplish this, we show that any sequence of these arrow maps would send $\omega c$ to $2ic$ before ever reaching $\omega'c$ for some positive odd integer $\omega'$.  After the first bottom map, we get $\omega c \overset{b}{\longrightarrow} 2kc + b - \omega c = \omega'c + b$ for some positive odd integer $\omega'$.  Certainly this is not an even multiple of $c$ since $b$ and $c$ share no common factors.  We now consider sequences of top maps and bottom maps, called \emph{double mappings}, by the meander.  This gives us two cases:


$$
\omega' c + xb \overset{ot}{\longrightarrow} 4kc - \omega'c - (x - 1)b \overset{b}{\longrightarrow} 2kc + xb - (4kc - \omega' c) = (\omega' - 2k) c + xb
$$
or
\beqs
\omega' c + xb & \overset{et}{\longrightarrow} & (2(2\ell - 1) - \omega')c - xb \\
~ & \overset{b}{\longrightarrow} & 2kc + b - ((2(2\ell - 1) - \omega')c - xb) \\
~ & = & (2k + \omega' - 2(2\ell - 1)) c + (x + 1)b
\eeqs
where $x$ is a positive integer.

In either case, the result is an odd multiple of $c$ plus a multiple of $b$ where the multiple of $b$ never decreases and increases by at most one after double mapping.  Since this double mapping can be repeated, we see that $\omega c$ will only map to $\omega' c + mb$ where $\omega'$ is odd and $m$ is a positive integer.  Suppose we have chosen $m$ to be the smallest positive integer with $\omega' c + mb = dc$ where $d$ is a positive integer.  Our goal is to show that $d$ is even.

Since $c$ and $b$ share no common factor, we know $b | (d - \omega')$ and $m = \frac{d - \omega'}{b} c$.  If $d$ is odd, then $d - \omega'$ is even, making $m$ even.  Then, if we consider $\frac{m}{2}$ as opposed to $m$, we see that 
$$
\omega' c + \frac{m}{2}b = dc - \frac{d - \omega'}{2} c = \frac{d - \omega'}{2} c
$$
which is an integer multiple of $c$, contradicting the minimality of the choice of $m$.  Thus, $d$ is even and every top-end-segment maps to the exterior face of the meander.  

Our next goal is to show that every segment maps to either a top-end-segment or the exterior face (in fact both are true).  This will complete the proof that there is no cycle in the meander and establish that the meander is Frobenius.

Consider an arbitrary segment and define an \emph{operation} to be first a bottom map and then a top map.  Note that since the total number of vertices is odd, there is no bottom-end-segment so an operation is well defined for all segments unless the bottom map arrives at a top-end-segment or the exterior face.  We will call such operations \emph{terminating} since the operation cannot actually be completed as defined.

As shown above, any operation where the top map does not occur within the last (odd) block must add $b$ plus an integer multiple of $c$ to the label of the segment.  Since we may suppose $b < 2a$ by Lemma~\ref{lemma:reduction} and we have assumed that there are at least $3$ blocks of size $a$ in this meander, every operation in which the top map occurs within the last block must be followed immediately by an operation in which the top map does not occur within the last block.  Thus, any sequence of operations must add an integer multiple of $c$ along with a strictly increasing positive integer multiple of $b$ to the label of the segment.  Since $b$ shares no common factors with $c$ (or $a$), this sequence of operations must terminate at either a top-end-segment or the exterior face, namely a segment label which is a multiple of $c$.  Thus, there can be no cycle in the meander so the meander is Frobenius.  \qed

\begin{figure}[ht!]
\begin{center}
 \epsfxsize 4in
 \epsfbox{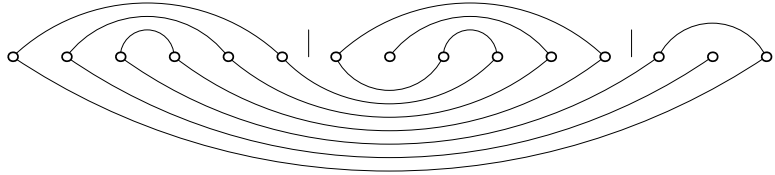}
 \caption{The meander $6|6|3$. \label{six-six-three}}
\end{center}
\end{figure}

\noindent \textit{Remark:} Although the first part of this proof may suggest that the result might hold for meanders of type $a|a|\dots |a|b$ where $a$ is even and $b$ is odd and possibly sharing a common factor with $a$, the example in Figure~\ref{six-six-three} is the smallest example that shows that this is not the case.  In this example, each of the top-end-segments do map to the exterior face but the other segments form a cycle since the operations never arrive at a multiple of $c$ and so never terminate.

Along with showing meanders are Frobenius, the more general problem of determining the index of an arbitrary meander is of some interest.  In particular, the results of this work lead to the following immediate observations.

\begin{fact}
A meander of type $2a|2a|\dots |2a|a$ has exactly $\ceil{\frac{a}{2}}$ components, of which exactly $\flr{\frac{a}{2}}$ are cycles. 
\end{fact}

\begin{fact}
A meander of type $a|a| \dots |a$ with a total of $n$ vertices has exactly $\ceil{\frac{n}{2a}} \ceil{\frac{a}{2}}$ components, of which exactly $\ceil{\frac{n}{2a}} \flr{\frac{a}{2}}$ are cycles.
\end{fact}

\section{Conclusion}

The classification of Frobenius Lie algebras appears to be a wild problem and the difficulty of classification may descend to seaweed Lie subalgebras of $\mathfrak{sl}_{n}$.  The wildness of the latter classification seems to present itself very quickly as the number of blocks increases.  Indeed, following Theorems~\ref{Thm:3-Blocks} and~\ref{Thm:4-Blocks}, one might expect that a Frobenius meander of type $a|b|c|d$ can be characterized by a relatively prime condition of the form
$$
\gcd (\alpha_{1} a + \alpha_{2} b + \alpha_{3} c + \alpha_{4} d, \beta_{1}a + \beta_{2}b + \beta_{3}c + \beta_{4}d) = 1,
$$
where the $\alpha_{i}$ and $\beta_{j}$ are integer \emph{coefficients}.  Substantial empirical evidence suggests that this is not so.  We argue as follows: Since the plane curve defined by a Frobenius meander is homotopically trivial, the meander can be contracted to a single point.  If such a ``winding down'' can be done through a deterministic sequence of edge contractions, each meander will yield a unique ``signature'', the elements of which represent the winding down moves performed.  This process can then be reversed, allowing us to build arbitrarily many Frobenius meanders with any given block structure.  The method of constructing such large sets of Frobenius meanders will be presented in a forthcoming paper by the current authors \cite{CMW12}.


Exhaustive simulations have shown that there is no set of coefficients, all with absolute value at most ten, that can be used in such a relatively prime condition.  All conditions with coefficients between negative ten and positive ten were checked against a set consisting of a large set of Frobenius meanders with five blocks.  No condition survived.  Since the addition of blocks seems to only complicate the situation, we are led to the following conjecture.

\begin{conjecture}\label{NoRelPrime}
There is no single relatively prime condition that suffices to classify Frobenius meanders with at least five blocks.
\end{conjecture}

Of course, if Conjecture~\ref{NoRelPrime} is true, it does not eliminate the possibility that some finite set of relatively prime conditions might classify meanders with five or more blocks that are Frobenius.  However, (less than exhaustive) simulations of numerous sets of conditions were tested against millions of Frobenius meanders.  No set survived.  We are led to the following stronger conjecture.

\begin{conjecture}
There is no finite set of relatively prime conditions that suffices to classify Frobenius meanders with at least five blocks.
\end{conjecture}


\end{document}